\documentclass[12pt]{article}
\usepackage{latexsym,amssymb,amsmath,ulem}
\usepackage{color}
\usepackage{graphicx}
\usepackage{latexsym,amssymb,amsmath,ulem}
\usepackage{amsfonts}
\usepackage{latexsym}
\usepackage{graphicx}
\usepackage{color}

\def\tto{\;{\lower 1pt \hbox{$\rightarrow$}}\kern -10pt
\hbox{\raise 2pt \hbox{$\rightarrow$}}\;}

\newtheorem{theorem}{Theorem}[section]
\newtheorem{proposition}{Proposition}[section]
\newtheorem{corollary}{Corollary}[section]
\newtheorem{lemma}{Lemma}[section]
\newtheorem{remark}{Remark}[section]
\newtheorem{example}{Example}[section]
\newtheorem{definition}{Definition}[section]

\numberwithin{equation}{section}

\renewcommand{\theequation}{\thesection.\arabic{equation}}
\normalsize

\makeatletter

\ifx\pdfoutput\relax\let\pdfoutput=\undefined\fi
\newcount\msipdfoutput
\ifx\pdfoutput\undefined
\else
 \ifcase\pdfoutput
 \else
    \ifx\paperwidth\undefined
    \else
      \ifdim\paperheight=0pt\relax
      \else
        \pdfpageheight\paperheight
      \fi
      \ifdim\paperwidth=0pt\relax
      \else
        \pdfpagewidth\paperwidth
      \fi
    \fi
  \fi
\fi

\newcount\@hour\newcount\@minute\chardef\@x10\chardef\@xv60
\def\tcitime{
\def\@time{%
  \@minute\time\@hour\@minute\divide\@hour\@xv
  \ifnum\@hour<\@x 0\fi\the\@hour:%
  \multiply\@hour\@xv\advance\@minute-\@hour
  \ifnum\@minute<\@x 0\fi\the\@minute
  }}%

\def\x@hyperref#1#2#3{%
   \catcode`\~ = 12
   \catcode`\$ = 12
   \catcode`\_ = 12
   \catcode`\# = 12
   \catcode`\& = 12
   \catcode`\% = 12
   \y@hyperref{#1}{#2}{#3}%
}

\def\y@hyperref#1#2#3#4{%
   #2\ref{#4}#3
   \catcode`\~ = 13
   \catcode`\$ = 3
   \catcode`\_ = 8
   \catcode`\# = 6
   \catcode`\& = 4
   \catcode`\% = 14
}

\@ifundefined{hyperref}{\let\hyperref\x@hyperref}{}
\@ifundefined{msihyperref}{\let\msihyperref\x@hyperref}{}

\@ifundefined{qExtProgCall}{\def\qExtProgCall#1#2#3#4#5#6{\relax}}{}
\def\QCTOpt[#1]#2{%
  \def\QCTOptB{#1}
  \def\QCTOptA{#2}
}
\def\QCTNOpt#1{%
  \def\QCTOptA{#1}
  \let\QCTOptB\empty
}
\def\Qct{%
  \@ifnextchar[{%
    \QCTOpt}{\QCTNOpt}
}
\def\QCBOpt[#1]#2{%
  \def\QCBOptB{#1}%
  \def\QCBOptA{#2}%
}
\def\QCBNOpt#1{%
  \def\QCBOptA{#1}%
  \let\QCBOptB\empty
}
\def\Qcb{%
  \@ifnextchar[{%
    \QCBOpt}{\QCBNOpt}%
}
\def\PrepCapArgs{%
  \ifx\QCBOptA\empty
    \ifx\QCTOptA\empty
      {}%
    \else
      \ifx\QCTOptB\empty
        {\QCTOptA}%
      \else
        [\QCTOptB]{\QCTOptA}%
      \fi
    \fi
  \else
    \ifx\QCBOptA\empty
      {}%
    \else
      \ifx\QCBOptB\empty
        {\QCBOptA}%
      \else
        [\QCBOptB]{\QCBOptA}%
      \fi
    \fi
  \fi
}
\newcount\GRAPHICSTYPE
\GRAPHICSTYPE=\z@
\def\GRAPHICSPS#1{%
 \ifcase\GRAPHICSTYPE
   \special{ps: #1}%
 \or
   \special{language "PS", include "#1"}%
 \fi
}%
\def\graffile#1#2#3#4{%
    \bgroup
	   \@inlabelfalse
       \leavevmode
       \@ifundefined{bbl@deactivate}{\def~{\string~}}{\activesoff}%
        \raise -#4 \BOXTHEFRAME{%
           \hbox to #2{\raise #3\hbox to #2{\null #1\hfil}}}%
    \egroup
}%
\def\draftbox#1#2#3#4{%
 \leavevmode\raise -#4 \hbox{%
  \frame{\rlap{\protect\tiny #1}\hbox to #2%
   {\vrule height#3 width\z@ depth\z@\hfil}%
  }%
 }%
}%
\newcount\@msidraft
\@msidraft=\z@
\let\nographics=\@msidraft
\newif\ifwasdraft
\wasdraftfalse

\def\GRAPHIC#1#2#3#4#5{%
   \ifnum\@msidraft=\@ne\draftbox{#2}{#3}{#4}{#5}%
   \else\graffile{#1}{#3}{#4}{#5}%
   \fi
}
\def\addtoLaTeXparams#1{%
    \edef\LaTeXparams{\LaTeXparams #1}}%
\newif\ifBoxFrame \BoxFramefalse
\newif\ifOverFrame \OverFramefalse
\newif\ifUnderFrame \UnderFramefalse

\def\BOXTHEFRAME#1{%
   \hbox{%
      \ifBoxFrame
         \frame{#1}%
      \else
         {#1}%
      \fi
   }%
}

\def\doFRAMEparams#1{\BoxFramefalse\OverFramefalse\UnderFramefalse\readFRAMEparams#1\end}%
\def\readFRAMEparams#1{%
 \ifx#1\end%
  \let\next=\relax
  \else
  \ifx#1i\dispkind=\z@\fi
  \ifx#1d\dispkind=\@ne\fi
  \ifx#1f\dispkind=\tw@\fi
  \ifx#1t\addtoLaTeXparams{t}\fi
  \ifx#1b\addtoLaTeXparams{b}\fi
  \ifx#1p\addtoLaTeXparams{p}\fi
  \ifx#1h\addtoLaTeXparams{h}\fi
  \ifx#1X\BoxFrametrue\fi
  \ifx#1O\OverFrametrue\fi
  \ifx#1U\UnderFrametrue\fi
  \ifx#1w
    \ifnum\@msidraft=1\wasdrafttrue\else\wasdraftfalse\fi
    \@msidraft=\@ne
  \fi
  \let\next=\readFRAMEparams
  \fi
 \next
 }%
\def\IFRAME#1#2#3#4#5#6{%
      \bgroup
      \let\QCTOptA\empty
      \let\QCTOptB\empty
      \let\QCBOptA\empty
      \let\QCBOptB\empty
      #6%
      \parindent=0pt
      \leftskip=0pt
      \rightskip=0pt
      \setbox0=\hbox{\QCBOptA}%
      \@tempdima=#1\relax
      \ifOverFrame
          \typeout{This is not implemented yet}%
          \show\HELP
      \else
         \ifdim\wd0>\@tempdima
            \advance\@tempdima by \@tempdima
            \ifdim\wd0 >\@tempdima
               \setbox1 =\vbox{%
                  \unskip\hbox to \@tempdima{\hfill\GRAPHIC{#5}{#4}{#1}{#2}{#3}\hfill}%
                  \unskip\hbox to \@tempdima{\parbox[b]{\@tempdima}{\QCBOptA}}%
               }%
               \wd1=\@tempdima
            \else
               \textwidth=\wd0
               \setbox1 =\vbox{%
                 \noindent\hbox to \wd0{\hfill\GRAPHIC{#5}{#4}{#1}{#2}{#3}\hfill}\\%
                 \noindent\hbox{\QCBOptA}%
               }%
               \wd1=\wd0
            \fi
         \else
            \ifdim\wd0>0pt
              \hsize=\@tempdima
              \setbox1=\vbox{%
                \unskip\GRAPHIC{#5}{#4}{#1}{#2}{0pt}%
                \break
                \unskip\hbox to \@tempdima{\hfill \QCBOptA\hfill}%
              }%
              \wd1=\@tempdima
           \else
              \hsize=\@tempdima
              \setbox1=\vbox{%
                \unskip\GRAPHIC{#5}{#4}{#1}{#2}{0pt}%
              }%
              \wd1=\@tempdima
           \fi
         \fi
         \@tempdimb=\ht1
         \advance\@tempdimb by -#2
         \advance\@tempdimb by #3
         \leavevmode
         \raise -\@tempdimb \hbox{\box1}%
      \fi
      \egroup%
}%
\def\DFRAME#1#2#3#4#5{%
  \vspace\topsep
  \hfil\break
  \bgroup
     \leftskip\@flushglue
	 \rightskip\@flushglue
	 \parindent\z@
	 \parfillskip\z@skip
     \let\QCTOptA\empty
     \let\QCTOptB\empty
     \let\QCBOptA\empty
     \let\QCBOptB\empty
	 \vbox\bgroup
        \ifOverFrame
           #5\QCTOptA\par
        \fi
        \GRAPHIC{#4}{#3}{#1}{#2}{\z@}%
        \ifUnderFrame
           \break#5\QCBOptA
        \fi
	 \egroup
  \egroup
  \vspace\topsep
  \break
}%
\def\FFRAME#1#2#3#4#5#6#7{%
  \@ifundefined{floatstyle}
    {
     \begin{figure}[#1]%
    }
    {
	 \ifx#1h
      \begin{figure}[H]%
	 \else
      \begin{figure}[#1]%
	 \fi
	}
  \let\QCTOptA\empty
  \let\QCTOptB\empty
  \let\QCBOptA\empty
  \let\QCBOptB\empty
  \ifOverFrame
    #4
    \ifx\QCTOptA\empty
    \else
      \ifx\QCTOptB\empty
        \caption{\QCTOptA}%
      \else
        \caption[\QCTOptB]{\QCTOptA}%
      \fi
    \fi
    \ifUnderFrame\else
      \label{#5}%
    \fi
  \else
    \UnderFrametrue%
  \fi
  \begin{center}\GRAPHIC{#7}{#6}{#2}{#3}{\z@}\end{center}%
  \ifUnderFrame
    #4
    \ifx\QCBOptA\empty
      \caption{}%
    \else
      \ifx\QCBOptB\empty
        \caption{\QCBOptA}%
      \else
        \caption[\QCBOptB]{\QCBOptA}%
      \fi
    \fi
    \label{#5}%
  \fi
  \end{figure}%
 }%
\newcount\dispkind%

\def\makeactives{
  \catcode`\"=\active
  \catcode`\;=\active
  \catcode`\:=\active
  \catcode`\'=\active
  \catcode`\~=\active
}
\bgroup
   \makeactives
   \gdef\activesoff{%
      \def"{\string"}%
      \def;{\string;}%
      \def:{\string:}%
      \def'{\string'}%
      \def~{\string~}%
    }
\egroup

\def\FRAME#1#2#3#4#5#6#7#8{%
 \bgroup
 \ifnum\@msidraft=\@ne
   \wasdrafttrue
 \else
   \wasdraftfalse%
 \fi
 \def\LaTeXparams{}%
 \dispkind=\z@
 \def\LaTeXparams{}%
 \doFRAMEparams{#1}%
 \ifnum\dispkind=\z@\IFRAME{#2}{#3}{#4}{#7}{#8}{#5}\else
  \ifnum\dispkind=\@ne\DFRAME{#2}{#3}{#7}{#8}{#5}\else
   \ifnum\dispkind=\tw@
    \edef\@tempa{\noexpand\FFRAME{\LaTeXparams}}%
    \@tempa{#2}{#3}{#5}{#6}{#7}{#8}%
    \fi
   \fi
  \fi
  \ifwasdraft\@msidraft=1\else\@msidraft=0\fi{}%
  \egroup
 }%
\def\TEXUX#1{"texux"}

\def\limfunc#1{\mathop{\rm #1}}%
\def\func#1{\mathop{\rm #1}\nolimits}%
\long\def\QQQ#1#2{%
     \long\expandafter\def\csname#1\endcsname{#2}}%
\@ifundefined{QTP}{\def\QTP#1{}}{}
\@ifundefined{QEXCLUDE}{\def\QEXCLUDE#1{}}{}
\@ifundefined{Qlb}{}{}
\@ifundefined{Qlt}{}{}
\long\def\QQA#1#2{}%
\def\QTR#1#2{{\csname#1\endcsname {#2}}}%

\def\EXPAND#1[#2]#3{}%
\def\NOEXPAND#1[#2]#3{}%
\def\LaTeXparent#1{}%
\def\ChildStyles#1{}%
\def\ChildDefaults#1{}%
\def\QTagDef#1#2#3{}%

\@ifundefined{correctchoice}{}{}
\@ifundefined{HTML}{\def\HTML#1{\relax}}{}
\@ifundefined{TCIIcon}{\def\TCIIcon#1#2#3#4{\relax}}{}
\if@compatibility
\else
  \providecommand{\UNICODE}[2][]{\protect\rule{.1in}{.1in}}
  \providecommand{\U}[1]{\protect\rule{.1in}{.1in}}
  
\fi

\@ifundefined{lambdabar}{
      
   }{}

\@ifundefined{StyleEditBeginDoc}{}{}
\def\QQfnmark#1{\footnotemark}

\@ifundefined{TCIMAKEINDEX}{}{\makeindex}%
\@ifundefined{abstract}{%
 \def\abstract{%
  \if@twocolumn
   \section*{Abstract (Not appropriate in this style!)}%
   \else \small
   \begin{center}{\bf Abstract\vspace{-.5em}\vspace{\z@}}\end{center}%
   \quotation
   \fi
  }%
 }{%
 }%
\@ifundefined{endabstract}{\def\endabstract
  {\if@twocolumn\else\endquotation\fi}}{}%
\@ifundefined{maketitle}{\def\maketitle#1{}}{}%
\@ifundefined{affiliation}{\def\affiliation#1{}}{}%
\@ifundefined{proof}{}{}%
\@ifundefined{endproof}{}{}%
\@ifundefined{newfield}{\def\newfield#1#2{}}{}%
\@ifundefined{chapter}{\def\chapter#1{\par(Chapter head:)#1\par }%
 \newcount\c@chapter}{}%
\@ifundefined{part}{\def\part#1{\par(Part head:)#1\par }}{}%
\@ifundefined{section}{\def\section#1{\par(Section head:)#1\par }}{}%
\@ifundefined{subsection}{\def\subsection#1%
 {\par(Subsection head:)#1\par }}{}%
\@ifundefined{subsubsection}{\def\subsubsection#1%
 {\par(Subsubsection head:)#1\par }}{}%
\@ifundefined{paragraph}{\def\paragraph#1%
 {\par(Subsubsubsection head:)#1\par }}{}%
\@ifundefined{subparagraph}{\def\subparagraph#1%
 {\par(Subsubsubsubsection head:)#1\par }}{}%
\@ifundefined{therefore}{}{}%
\@ifundefined{backepsilon}{}{}%
\@ifundefined{yen}{}{}%
\@ifundefined{registered}{%
   \def\registered{\relax\ifmmode{}\r@gistered
                    \else$\m@th\r@gistered$\fi}%
 \def\r@gistered{^{\ooalign
  {\hfil\raise.07ex\hbox{$\scriptstyle\rm\text{R}$}\hfil\crcr
  \mathhexbox20D}}}}{}%
\@ifundefined{Eth}{}{}%
\@ifundefined{eth}{}{}%
\@ifundefined{Thorn}{}{}%
\@ifundefined{thorn}{}{}%
\@ifundefined{degree}{}{}%
\newdimen\theight
\@ifundefined{Column}{\def\Column{%
 \vadjust{\setbox\z@=\hbox{\scriptsize\quad\quad tcol}%
  \theight=\ht\z@\advance\theight by \dp\z@\advance\theight by \lineskip
  \kern -\theight \vbox to \theight{%
   \rightline{\rlap{\box\z@}}%
   \vss
   }%
  }%
 }}{}%
\@ifundefined{qed}{\def\qed{%
 \ifhmode\unskip\nobreak\fi\ifmmode\ifinner\else\hskip5\p@\fi\fi
 \hbox{\hskip5\p@\vrule width4\p@ height6\p@ depth1.5\p@\hskip\p@}%
 }}{}%
\@ifundefined{cents}{}{}%
\@ifundefined{tciLaplace}{}{}%
\@ifundefined{tciFourier}{}{}%
\@ifundefined{textcurrency}{}{}%
\@ifundefined{texteuro}{}{}%
\@ifundefined{euro}{}{}%
\@ifundefined{textfranc}{}{}%
\@ifundefined{textlira}{}{}%
\@ifundefined{textpeseta}{}{}%
\@ifundefined{miss}{\def\miss{\hbox{\vrule height2\p@ width 2\p@ depth\z@}}}{}%
\@ifundefined{vvert}{}{}
\@ifundefined{tcol}{\def\tcol#1{{\baselineskip=6\p@ \vcenter{#1}} \Column}}{}%
\@ifundefined{dB}{}{}
\@ifundefined{mB}{}{}
\@ifundefined{nB}{}{}
\@ifundefined{note}{}{}%
\def\newfmtname{LaTeX2e}
\ifx\fmtname\newfmtname
  \DeclareOldFontCommand{\rm}{\normalfont\rmfamily}{\mathrm}
  \DeclareOldFontCommand{\sf}{\normalfont\sffamily}{\mathsf}
  \DeclareOldFontCommand{\tt}{\normalfont\ttfamily}{\mathtt}
  \DeclareOldFontCommand{\bf}{\normalfont\bfseries}{\mathbf}
  \DeclareOldFontCommand{\it}{\normalfont\itshape}{\mathit}
  \DeclareOldFontCommand{\sl}{\normalfont\slshape}{\@nomath\sl}
  \DeclareOldFontCommand{\sc}{\normalfont\scshape}{\@nomath\sc}
\fi

\def\alpha{{\Greekmath 010B}}%
\def\beta{{\Greekmath 010C}}%
\def\gamma{{\Greekmath 010D}}%
\def\delta{{\Greekmath 010E}}%
\def\epsilon{{\Greekmath 010F}}%
\def\zeta{{\Greekmath 0110}}%
\def\eta{{\Greekmath 0111}}%
\def\theta{{\Greekmath 0112}}%
\def\iota{{\Greekmath 0113}}%
\def\kappa{{\Greekmath 0114}}%
\def\lambda{{\Greekmath 0115}}%
\def\mu{{\Greekmath 0116}}%
\def\nu{{\Greekmath 0117}}%
\def\xi{{\Greekmath 0118}}%
\def\pi{{\Greekmath 0119}}%
\def\rho{{\Greekmath 011A}}%
\def\sigma{{\Greekmath 011B}}%
\def\tau{{\Greekmath 011C}}%
\def\upsilon{{\Greekmath 011D}}%
\def\phi{{\Greekmath 011E}}%
\def\chi{{\Greekmath 011F}}%
\def\psi{{\Greekmath 0120}}%
\def\omega{{\Greekmath 0121}}%
\def\varepsilon{{\Greekmath 0122}}%
\def\vartheta{{\Greekmath 0123}}%
\def\varpi{{\Greekmath 0124}}%
\def\varrho{{\Greekmath 0125}}%
\def\varsigma{{\Greekmath 0126}}%
\def\varphi{{\Greekmath 0127}}%

\def\nabla{{\Greekmath 0272}}
\def\FindBoldGroup{%
   {\setbox0=\hbox{$\mathbf{x\global\edef\theboldgroup{\the\mathgroup}}$}}%
}

\def\Greekmath#1#2#3#4{%
    \if@compatibility
        \ifnum\mathgroup=\symbold
           \mathchoice{\mbox{\boldmath$\displaystyle\mathchar"#1#2#3#4$}}%
                      {\mbox{\boldmath$\textstyle\mathchar"#1#2#3#4$}}%
                      {\mbox{\boldmath$\scriptstyle\mathchar"#1#2#3#4$}}%
                      {\mbox{\boldmath$\scriptscriptstyle\mathchar"#1#2#3#4$}}%
        \else
           \mathchar"#1#2#3#4%
        \fi
    \else
        \FindBoldGroup
        \ifnum\mathgroup=\theboldgroup 
           \mathchoice{\mbox{\boldmath$\displaystyle\mathchar"#1#2#3#4$}}%
                      {\mbox{\boldmath$\textstyle\mathchar"#1#2#3#4$}}%
                      {\mbox{\boldmath$\scriptstyle\mathchar"#1#2#3#4$}}%
                      {\mbox{\boldmath$\scriptscriptstyle\mathchar"#1#2#3#4$}}%
        \else
           \mathchar"#1#2#3#4%
        \fi     	
	  \fi}

\newif\ifGreekBold  \GreekBoldfalse
\let\SAVEPBF=\pbf
\def\pbf{\GreekBoldtrue\SAVEPBF}%

\@ifundefined{theorem}{\newtheorem{theorem}{Theorem}}{}
\@ifundefined{lemma}{\newtheorem{lemma}[theorem]{Lemma}}{}
\@ifundefined{corollary}{}{}
\@ifundefined{conjecture}{}{}
\@ifundefined{proposition}{\newtheorem{proposition}[theorem]{Proposition}}{}
\@ifundefined{axiom}{}{}
\@ifundefined{remark}{\newtheorem{remark}{Remark}}{}
\@ifundefined{example}{\newtheorem{example}{Example}}{}
\@ifundefined{exercise}{}{}
\@ifundefined{definition}{}{}

\@ifundefined{mathletters}{%
  \newcounter{equationnumber}
  \def\mathletters{%
     \addtocounter{equation}{1}
     \edef\@currentlabel{\theequation}%
     \setcounter{equationnumber}{\c@equation}
     \setcounter{equation}{0}%
     \edef\theequation{\@currentlabel\noexpand\alph{equation}}%
  }
  
}{}

\@ifundefined{BibTeX}{%
    \def\BibTeX{{\rm B\kern-.05em{\sc i\kern-.025em b}\kern-.08em
                 T\kern-.1667em\lower.7ex\hbox{E}\kern-.125emX}}}{}%
\@ifundefined{AmS}%
    {\def\AmS{{\protect\usefont{OMS}{cmsy}{m}{n}%
                A\kern-.1667em\lower.5ex\hbox{M}\kern-.125emS}}}{}%
\@ifundefined{AmSTeX}{}{}%
\def\@@eqncr{\let\@tempa\relax
    \ifcase\@eqcnt \def\@tempa{& & &}\or \def\@tempa{& &}%
      \else \def\@tempa{&}\fi
     \@tempa
     \if@eqnsw
        \iftag@
           \@taggnum
        \else
           \@eqnnum\stepcounter{equation}%
        \fi
     \fi
     \global\tag@false
     \global\@eqnswtrue
     \global\@eqcnt\z@\cr}

\def\TCItag{\@ifnextchar*{\@TCItagstar}{\@TCItag}}
\def\@TCItag#1{%
    \global\tag@true
    \global\def\@taggnum{(#1)}%
    \global\def\@currentlabel{#1}}
\def\@TCItagstar*#1{%
    \global\tag@true
    \global\def\@taggnum{#1}%
    \global\def\@currentlabel{#1}}
\def\tint{\msi@int\textstyle\int}%
\def\tiint{\msi@int\textstyle\iint}%
\def\tiiint{\msi@int\textstyle\iiint}%
\def\tiiiint{\msi@int\textstyle\iiiint}%
\def\tidotsint{\msi@int\textstyle\idotsint}%
\def\toint{\msi@int\textstyle\oint}%

\def\tbigcap{\mathop{\textstyle \bigcap }}%
\newtoks\temptoksa
\newtoks\temptoksb
\newtoks\temptoksc

\def\msi@int#1#2{%
 \def\@temp{{#1#2\the\temptoksc_{\the\temptoksa}^{\the\temptoksb}}}%
 \futurelet\@nextcs
 \@int
}

\def\@int{%
   \ifx\@nextcs\limits
      \typeout{Found limits}%
      \temptoksc={\limits}%
	  \let\@next\@intgobble%
   \else\ifx\@nextcs\nolimits
      \typeout{Found nolimits}%
      \temptoksc={\nolimits}%
	  \let\@next\@intgobble%
   \else
      \typeout{Did not find limits or no limits}%
      \temptoksc={}%
      \let\@next\msi@limits%
   \fi\fi
   \@next
}%

\def\@intgobble#1{%
   \typeout{arg is #1}%
   \msi@limits
}

\def\msi@limits{%
   \temptoksa={}%
   \temptoksb={}%
   \@ifnextchar_{\@limitsa}{\@limitsb}%
}

\def\@limitsa_#1{%
   \temptoksa={#1}%
   \@ifnextchar^{\@limitsc}{\@temp}%
}

\def\@limitsb{%
   \@ifnextchar^{\@limitsc}{\@temp}%
}

\def\@limitsc^#1{%
   \temptoksb={#1}%
   \@ifnextchar_{\@limitsd}{\@temp}%
}

\def\@limitsd_#1{%
   \temptoksa={#1}%
   \@temp
}

\def\dint{\msi@int\displaystyle\int}%
\def\diint{\msi@int\displaystyle\iint}%
\def\diiint{\msi@int\displaystyle\iiint}%
\def\diiiint{\msi@int\displaystyle\iiiint}%
\def\didotsint{\msi@int\displaystyle\idotsint}%
\def\doint{\msi@int\displaystyle\oint}%

\if@compatibility\else
  \RequirePackage{amsmath}
\fi

\begin{document}

\title{Duality for convex infinite optimization on linear spaces}
\author{M. A. Goberna\thanks{%
Department of Mathematics, University of Alicante, Alicante, Spain
(mgoberna@ua.es). Corresponding author.} \ and\ M. Volle\thanks{%
Avignon University, LMA EA 2151, Avignon, France
(michel.volle@univ-avignon.fr)} }
\maketitle
\date{}

\begin{abstract}
This note establishes a limiting formula for the conic Lagrangian dual of a
convex infinite optimization problem, correcting the classical version of
Karney [Math. Programming 27 (1983) 75-82] for convex semi-infinite
programs. A reformulation of the convex infinite optimization problem with a
single constraint leads to a limiting formula for the corresponding
Lagrangian dual, called sup-dual, and also for the primal problem in the
case when strong Slater condition holds, which also entails strong
sup-duality.
\end{abstract}




\section{Introduction}

Given a real linear space $X,$ consider the (algebraic) convex infinite
programming (CIP) problem
\begin{equation*}
(P)\text{ }\inf_{x\in X}f(x),\text{ }\mathrm{s.t.}\text{ }f_{t}(x)\leq 0,\
t\in T,
\end{equation*}%
where $T$ is an infinite index set and $f,f_{t}:X\longrightarrow \overline{%
\mathbb{R}}:=\mathbb{R\cup }\left\{ \pm \infty \right\} ,$ $t\in T,$\ are
convex proper functions. We denote by%
\begin{equation*}
E:=\tbigcap\limits_{t\in T}\ \left[ f_{t}\leq 0\right] =\left\{ x\in
X:f_{t}(x)\leq 0,\ t\in T\right\}
\end{equation*}%
the feasible set of $(P)$ and define
\begin{equation*}
M:=\tbigcap\limits_{t\in T}\ \limfunc{dom}f_{t}\supset E\text{ and }\Delta
:=M\cap \limfunc{dom}f.
\end{equation*}

Let $\mathbb{R}_{+}^{(T)}$ be the positive cone of the space $\mathbb{R}%
^{(T)}$ of functions\textbf{\ }$\lambda =\left( \lambda \right) _{t\in
T}:T\rightarrow \mathbb{R}$ whose support $\limfunc{supp}\lambda :=\left\{
t\in T:\lambda _{t}\neq 0\right\} $ is finite and let $0_{\mathbb{R}^{(T)}}$
be its null element. The ordinary \textit{Lagrangian} \textit{function}
associated to $(P)$ is (see \cite{GL98}, \cite{Karney83}, etc.) is $%
L_{0}:X\times \mathbb{R}_{+}^{(T)}\longrightarrow \overline{\mathbb{R}}$
such that $L_{0}\left( x,\lambda \right) :=f(x)+\sum_{t\in T}\lambda
_{t}f_{t}(x),$ where

\begin{equation*}
\sum_{t\in T}\lambda _{t}f_{t}(x):=\left\{
\begin{array}{ll}
\sum_{t\in \limfunc{supp}\lambda }\lambda _{t}f_{t}(x), & \text{if }\lambda
\neq 0_{\mathbb{R}^{(T)}}, \\
0, & \text{if }\lambda =0_{\mathbb{R}^{(T)}}.%
\end{array}%
\right.
\end{equation*}

A slightly different Lagrangian is the associated to the cone constrained
reformulation of $(P),$\ that is \cite[page 138]{Za02}, the function $%
L:X\times \mathbb{R}_{+}^{(T)}\longrightarrow \overline{\mathbb{R}}$ such
that%
\begin{equation*}
L\left( x,\lambda \right) :=\left\{
\begin{array}{ll}
f(x)+\sum_{t\in T}\lambda _{t}f_{t}(x), & \text{if }x\in M,\text{ }\lambda
\in \mathbb{R}_{+}^{(T)}, \\
+\infty , & \text{else.}%
\end{array}%
\right.
\end{equation*}%
We call $L$ the \textit{conic Lagrangian} of $(P).$

For each $x\in X$ we have
\begin{equation*}
\sup_{\lambda \in \mathbb{R}_{+}^{(T)}}L_{0}\left( x,\lambda \right)
=\sup_{\lambda \in \mathbb{R}_{+}^{(T)}}L\left( x,\lambda \right)
=f(x)+\delta _{E}(x),
\end{equation*}%
where $\delta _{E}$ is the indicator of $E,$ that is, $\delta _{E}\left(
x\right) =0$ if $x\in E$ and $\delta _{E}\left( x\right) =+\infty $
otherwise. Consequently,%
\begin{equation*}
\inf_{x\in X}\sup_{\lambda \in \mathbb{R}_{+}^{(T)}}L_{0}\left( x,\lambda
\right) =\inf_{x\in X}\sup_{\lambda \in \mathbb{R}_{+}^{(T)}}L\left(
x,\lambda \right) =\inf (P).
\end{equation*}%
The \textit{ordinary} and \textit{conic-Lagrangian dual problems} of $(P)$
read, respectively,
\begin{equation*}
(D_{0})\text{ }\sup_{\lambda \in \mathbb{R}_{+}^{(T)}}\inf_{x\in X}\left(
f(x)+\sum_{t\in T}\lambda _{t}f_{t}(x)\right) ,\text{ }
\end{equation*}%
and%
\begin{equation*}
(D)\text{ }\sup_{\lambda \in \mathbb{R}_{+}^{(T)}}\inf_{x\in M}\left(
f(x)+\sum_{t\in T}\lambda _{t}f_{t}(x)\right) ,
\end{equation*}%
and one has%
\begin{equation}
\sup (D_{0})\leq \sup (D)\leq \inf \left( P\right) .  \label{1.1}
\end{equation}

Note that, if $\limfunc{dom}f\subset M,$ then $\sup (D_{0})=\sup (D).$ This
is in particular the case when the functions $f_{t},$ $t\in T,$ are
real-valued. But it may happen that $\sup (D_{0})<\sup (D)$ even if $T$ is
finite and Slater condition holds. This is the case in the next example.

\begin{example}
\label{Exam1} Consider $X=\mathbb{R}^{2},$ $T=\left\{ 1\right\} ,$ $f\left(
x_{1},x_{2}\right) =e^{x_{2}},$ and
\begin{equation*}
f_{1}\left( x_{1},x_{2}\right) =\left\{
\begin{array}{ll}
x_{1}, & \text{if }x_{2}\geq 0, \\
+\infty , & \text{if }x_{2}<0.%
\end{array}%
\right.
\end{equation*}%
We then have%
\begin{equation*}
\max (D_{0})=0<1=\max (D)=\min \left( P\right) .
\end{equation*}
\end{example}

Duffin \cite{Duffin73} observed that a positive duality gap may occur when
one considers the ordinary Lagrangian dual $(D_{0})$ of $\left( P\right) .$
The same happens when $(D_{0})$ is replaced by $(D)$ even though, according
to (\ref{1.1}), the gap may be smaller.\ Different ways have been proposed
to close the duality gap, e.g., by adding a linear perturbation to the
saddle function $f+\sum_{t\in T}\lambda _{t}f_{t},$ and sending it to zero
in the limit \cite{Duffin73}. Blair, Duffin and Jeroslow \cite{BDJ82} used
the conjugate duality theory to extend the limiting phenomena to the general
minimax setting. Pomerol \cite{Pomerol82} showed that it was possible to
obtain infisup theorems, including that of \cite{BDJ82}, by using a slightly
more general form of the duality theory. In turn, Karney and Morley \cite%
{KM86} proved that, when $X=\mathbb{R}^{n},$ either the convex semi-infinite
programming (CSIP in brief) problem $\left( P\right) $ satisfies some
recession condition guaranteeing a zero duality gap or there exists $d\in
\mathbb{R}^{n}\diagdown \left\{ 0_{n}\right\} $ such that the problem
\begin{equation*}
\left( P_{\varepsilon }\right) \text{ }\inf_{x\in X}f(x)+\varepsilon
\left\langle d,x\right\rangle ,\text{ }\mathrm{s.t.}\text{ }f_{t}(x)\leq 0,\
t\in T,
\end{equation*}%
satisfies the mentioned recession condition for $\varepsilon >0$
sufficiently small, with $\left( P_{\varepsilon }\right) $ enjoying strong
duality, and $\inf \left( P\right) =\lim\limits_{\varepsilon \downarrow
0}\left( P_{\varepsilon }\right) .$ The theory developed in \cite{KM86}
subsumed the CSIP versions of some results on limiting Lagrangians in \cite%
{Borwein80} and \cite{DJ79}. Three years before, Karney gave, in the CSIP
setting, a limiting formula for the dual problem $(D_{0}):$%
\begin{equation}
\sup (D_{0})=\lim\limits_{\varepsilon \downarrow 0}\inf \left\{ f(x):\text{ }%
f_{t}(x)\leq \varepsilon ,\ t\in T\right\} .  \label{1.2}
\end{equation}%
According to \cite[Proposition 3.1]{Karney83}, this formula comes from \cite[%
Theorem 7]{Rockafellar74} and \cite[Corollary 2]{Borwein80}, and does not
require any constraint qualification (other than $E\neq \emptyset ,$ or
something stronger as $E\cap \func{dom}f\neq \emptyset ,$ $E\subset \limfunc{%
cl}\func{dom}f,$ ...). The next example shows that \cite[Proposition 3.1]%
{Karney83} fails even in linear semi-infinite programming, where $\func{dom}%
f=X=\mathbb{R}^{n}.$

\begin{example}
\label{Exam2}Consider the following optimization problem, with $T=\mathbb{N}%
: $
\begin{equation*}
\begin{array}{llll}
\left( P\right) & \inf_{x\in \mathbb{R}^{2}} & x_{2}\qquad &  \\
& \text{\textrm{s.t.}} & x_{1}\leq 0, & \left( t=1\right) \\
&  & -x_{2}\leq 1, & \left( t=2\right) \\
&  & t^{-1}x_{1}-x_{2}\leq 0, & t=3,4,...%
\end{array}%
\end{equation*}%
Its dual problem $\left( D_{0}\right) $ is equivalent to the Haar dual (see,
e.g., \cite{GL98})%
\begin{equation*}
\begin{array}{ll}
\sup_{\lambda \in \mathbb{R}_{+}^{(\mathbb{N})}} & -\lambda _{2} \\
\text{\textrm{s.t.}} & \lambda _{1}\left(
\begin{array}{c}
-1 \\
0%
\end{array}%
\right) +\lambda _{2}\left(
\begin{array}{c}
0 \\
1%
\end{array}%
\right) +\sum\nolimits_{t\geq 3}\left(
\begin{array}{c}
-t^{-1} \\
1%
\end{array}%
\right) =\left(
\begin{array}{c}
0 \\
1%
\end{array}%
\right) ,%
\end{array}%
\end{equation*}%
whose unique\ feasible solution is $\lambda \in \mathbb{R}_{+}^{(\mathbb{N}%
)} $ such that $\lambda _{2}=1$ and $\lambda _{t}=0$ for $t\neq 2.$ So, $%
\max \left( D_{0}\right) =-1$ while $E=\left\{ \left( x_{1},x_{2}\right)
:x_{1}\leq 0,x_{2}\geq 0\right\} ,$ so that $\min \left( P\right) =0.$ On
the other hand, given $\varepsilon >0,$
\begin{equation*}
\left\{ x\in \mathbb{R}^{2}:\text{ }f_{t}(x)\leq \varepsilon ,\ t\in \mathbb{%
N}\right\} =\left\{ x\in \mathbb{R}^{2}:x_{1}\leq \varepsilon ,x_{2}\geq
-\varepsilon ,\frac{x_{1}}{3}-x_{2}\leq \varepsilon \right\} ,
\end{equation*}%
so that
\begin{equation*}
\min \left\{ x_{2}:f_{t}(x)\leq \varepsilon ,\ t\in \mathbb{N}\right\}
=-\varepsilon
\end{equation*}%
is attained at $\left\{ \left( x_{1},-\varepsilon \right) :x_{1}\leq
0\right\} .$ Hence,%
\begin{equation*}
\max \left( D_{0}\right) =-1<0=\lim\limits_{\varepsilon \downarrow 0}\min
\left\{ x_{2}:f_{t}(x)\leq \varepsilon ,\ t\in \mathbb{N}\right\} .
\end{equation*}
\end{example}

From \cite[Proposition 3.1]{Karney83} Karney obtained, following the
suggestion of an unknown referee, the reverse stromg duality theorem \cite[%
Theorem 3.2]{Karney83}%
\begin{equation*}
\min \left( P\right) =\sup \left( D_{0}\right)
\end{equation*}%
under some recession condition. However, he asserted in \cite[Section 5]%
{Karney83} that he had two (longer) unpublished proofs. In either case, his
result has been recently proved from a new strong duality theorem for CIP
(see \cite[Corollary 3.2 and Remark 3.2]{DGLV21}).

In this note we show in a simpler way, for general CIP problems, that, if
\begin{equation*}
\exists \alpha >0,\exists a\in \func{dom}f:\ \ f_{t}(a)\leq -\alpha ,\
\forall t\in T,
\end{equation*}%
then (\ref{1.2}) entails that zero duality gap holds:%
\begin{equation*}
\sup (D_{0})=\inf \left( P\right) .
\end{equation*}%
This duality theorem is obtained by studying the Lagrangian dual $(D_{1})$
associated with the representation of $E$ by a single constraint (the
so-called sup-function). Section 2 (resp. Section 3) provides a limiting
formula for $\sup (D)$ (resp. $\sup (D_{1})$). Under the strong Slater
condition, the limiting formula for $\sup (D_{1})$ also holds for $\inf
\left( P\right) $ together with the strong duality theorem $\inf \left(
P\right) =\max (D_{1}).$

\section{Conic-Lagrangian duality}

Problem $(D)$ receives a perturbational interpretation (see \cite{Bot10},
\cite{Za02}, etc.) in terms of the \textit{ordinary value function} $v:%
\mathbb{R}^{T}\longrightarrow \overline{\mathbb{R}}$ associated with $(P)$
defined by%
\begin{equation*}
v\left( y\right) :=\inf \left\{ f(x):f_{t}(x)\leq y_{t},t\in T\right\}
,\forall y=\left( y_{t}\right) _{t\in T}\in \mathbb{R}^{T}.
\end{equation*}

Let us make explicit this approach. The linear space $Y:=\mathbb{R}^{T},$
equipped with the product topology, is a locally convex Hausdorff
topological vector space whose topological dual is $\mathbb{R}^{(T)}$ via
the bilinear pairing%
\begin{equation*}
\left\langle \cdot ,\cdot \right\rangle :Y\times \mathbb{R}%
^{(T)}\longrightarrow \mathbb{R}\text{ such that }\left\langle y,\lambda
\right\rangle =\sum_{t\in T}\lambda _{t}y_{t}.
\end{equation*}

The Fenchel conjugate of $v$ is (see \cite{Bot10}, \cite{Za02}, etc.)
\begin{equation}
-v^{\ast }\left( -\lambda \right) =\left\{
\begin{array}{ll}
\inf_{x\in \Delta }\left( f\left( x\right) +\sum_{t\in T}\lambda
_{t}f_{t}\left( x\right) \right) , & \text{if }\Delta \neq \emptyset \text{
and }\lambda \in \mathbb{R}_{+}^{(T)}, \\
-\infty , & \text{if }\Delta =\emptyset \text{ or }\lambda \in \mathbb{R}%
^{(T)}\diagdown \mathbb{R}_{+}^{(T)}.%
\end{array}%
\right.  \label{2.1}
\end{equation}%
If $\Delta \neq \emptyset $ we the have
\begin{equation*}
\begin{array}{ll}
v^{\ast \ast }\left( 0_{Y}\right) & =\sup_{\lambda \in \mathbb{R}%
^{(T)}}-v^{\ast }\left( \lambda \right) =\sup_{\lambda \in \mathbb{R}%
^{(T)}}-v^{\ast }\left( -\lambda \right) =\sup_{\lambda \in \mathbb{R}%
_{+}^{(T)}}-v^{\ast }\left( -\lambda \right) \\
& =\sup_{\lambda \in \mathbb{R}_{+}^{(T)}}\inf_{x\in \Delta }\left( f\left(
x\right) +\sum_{t\in T}\lambda _{t}f_{t}\left( x\right) \right) =\sup (D).%
\end{array}%
\end{equation*}%
Note that, if $\Delta =\emptyset $ we have $\limfunc{dom}v=\emptyset $ and $%
v^{\ast \ast }\left( 0_{Y}\right) =+\infty =\sup (D).$ Therefore, in all
cases we have%
\begin{equation}
\sup (D)=v^{\ast \ast }\left( 0_{Y}\right) \leq \overline{v}\left(
0_{Y}\right) \leq v\left( 0_{Y}\right) =\inf \left( P\right) ,  \label{2.3}
\end{equation}%
where $\overline{v}$ is the lower semicontinuous (lsc in brief) hull of $v$
for the product topology on $Y=\mathbb{R}^{T}.$ A neighborhood basis of the
origin $0_{Y}$ is furnished by the family
\begin{equation*}
\left\{ V_{\varepsilon }^{H}:\varepsilon >0,H\in \mathcal{F}\left( T\right)
\right\} ,
\end{equation*}%
where $\mathcal{F}\left( T\right) $ is the class of non-empty finite subsets
of $T,$ and
\begin{equation*}
V_{\varepsilon }^{H}:=\left\{ y\in Y:\left\vert y_{t}\right\vert \leq
\varepsilon ,t\in H\right\} .
\end{equation*}%
We now give a general explicit formula for $\overline{v}\left( 0_{Y}\right)
: $

\begin{lemma}
\label{lemma2.1} $\overline{v}\left( 0_{Y}\right) =\sup\limits_{\varepsilon
>0,H\in \mathcal{F}\left( T\right) }\inf\limits_{x\in M}\left\{ f\left(
x\right) :f_{t}\left( x\right) \leq \varepsilon ,t\in H\right\} .$
\end{lemma}

\noindent \textbf{Proof} For each $\varepsilon >0$ and $H\in \mathcal{F}%
\left( T\right) $ one has%
\begin{equation*}
\begin{array}{ll}
\inf\limits_{y\in V_{\varepsilon }^{H}}v\left( y\right) & =\inf \left\{
f\left( x\right) :f_{t}\left( x\right) \leq y_{t},t\in T;\left\vert
y_{t}\right\vert \leq \varepsilon ,t\in H\right\} \\
& =\inf \left\{ f\left( x\right) :f_{t}\left( x\right) \leq \varepsilon
,t\in H;f_{t}\left( x\right) <+\infty ,t\notin H\right\} \\
& =\inf\limits_{x\in M}\left\{ f\left( x\right) :f_{t}\left( x\right) \leq
\varepsilon ,t\in H\right\} .%
\end{array}%
\end{equation*}%
Since $\overline{v}\left( 0_{Y}\right) =\sup\limits_{\varepsilon >0,H\in
\mathcal{F}\left( T\right) }\inf\limits_{y\in V_{\varepsilon }^{H}}v\left(
y\right) ,$ we are done.\hfill $\square \medskip $

\begin{remark}
\label{Rem 2.1} From Lemma \ref{lemma2.1} one gets%
\begin{equation*}
\overline{v}\left( 0_{Y}\right) \leq \lim\limits_{\varepsilon \downarrow
0}\inf \left\{ f\left( x\right) :f_{t}\left( x\right) \leq \varepsilon ,t\in
T\right\} .
\end{equation*}
\end{remark}

\hfill

\begin{remark}
\label{Rem2.2} In the case when the index set $T$ is finite, the formula
provided by Lemma \ref{lemma2.1} can be simplified as follows:%
\begin{equation*}
\overline{v}\left( 0_{Y}\right) =\lim\limits_{\varepsilon \downarrow 0}\inf
\left\{ f\left( x\right) :f_{t}\left( x\right) \leq \varepsilon ,t\in
T\right\} .
\end{equation*}%
In such a case we also have $M=\tbigcap\limits_{t\in T}\ \limfunc{dom}f_{t}$
and%
\begin{equation*}
v^{\ast \ast }\left( 0_{Y}\right) =\sup_{\lambda \in \mathbb{R}%
_{+}^{T}}\inf_{x\in M}\left( f\left( x\right) +\sum_{t\in T}\lambda
_{t}f_{t}\left( x\right) \right) .
\end{equation*}
\end{remark}

\begin{proposition}[Limiting formula for $\sup (D)$]
\label{Prop2.1} Assume either $\overline{v}\left( 0_{Y}\right) \neq +\infty $
or $\sup (D)\neq -\infty .$ Then we have%
\begin{equation*}
\sup (D)=\sup\limits_{\varepsilon >0,H\in \mathcal{F}\left( T\right)
}\inf\limits_{x\in M}\left\{ f\left( x\right) :f_{t}\left( x\right) \leq
\varepsilon ,t\in H\right\} .
\end{equation*}
\end{proposition}

\noindent \textbf{Proof} We know that $\sup (D)=v^{\ast \ast }\left(
0_{Y}\right) $ (see (\ref{2.3})). Since the functions $f$ and $f_{t},$ $t\in
T,$ are convex, the value function $v$ is convex, too. By \cite[Proposition 1%
]{Borwein80}, we then have $\sup (D)=\overline{v}\left( 0_{Y}\right) $ and
Lemma \ref{lemma2.1} concludes the proof.\hfill $\square \medskip $

\begin{remark}
\label{Rem2.3} Condition $\overline{v}\left( 0_{Y}\right) \neq +\infty $ is
in particular satisfied if $\inf (P)\neq +\infty ,$ that is $E\cap \limfunc{%
dom}f\neq \emptyset .$ \newline
Condition $\sup (D)\neq -\infty $ is satisfied if and only if there exists $%
\lambda \in \mathbb{R}_{+}^{(T)}$ and $r\in \mathbb{R}$ such that%
\begin{equation*}
x\in M\Longrightarrow f\left( x\right) +\sum_{t\in T}\lambda _{t}f_{t}\left(
x\right) \geq r.
\end{equation*}
\end{remark}

\begin{remark}
\label{Rem2.4} By (\ref{1.1}), (\ref{2.1}) and (\ref{2.3}), we have%
\begin{equation*}
\sup (D_{0})\leq \sup (D)\leq \lim\limits_{\varepsilon \downarrow 0}\inf
\left\{ f\left( x\right) :f_{t}\left( x\right) \leq \varepsilon ,t\in
T\right\} .
\end{equation*}%
In \cite[Proposition 3.1]{Karney83} it is claimed that for $X=\mathbb{R}%
^{n}, $ $f$ and $f_{t},$ $t\in T,$ are proper, lsc and convex, and $E\neq
\emptyset ,$ it holds that
\begin{equation*}
\sup (D_{0})=\lim\limits_{\varepsilon \downarrow 0}\inf \left\{ f\left(
x\right) :f_{t}\left( x\right) \leq \varepsilon ,t\in T\right\} .
\end{equation*}%
To the best of our knowledge, this fact has not been proved anywhere. We
prove in Proposition \ref{Prop3.2} below an exact formula for its right-hand
side.
\end{remark}

\section{Sup-Lagrangian duality}

Let $h:=\sup\limits_{t\in T}f_{t}$ be the \textit{sup-function} of $\left(
P\right) $ which allows to represent its feasible set $E$ with a single
constraint. We associate with $\left( P\right) $ another Lagrangian $%
L_{1}:X\times \mathbb{R}_{+}\longrightarrow \overline{\mathbb{R}},$ called
\textit{sup-Lagrangian,} such that%
\begin{equation*}
L_{1}\left( x,s\right) :=\left\{
\begin{array}{ll}
f(x)+sh(x), & \text{if }x\in \Delta _{1}:=\limfunc{dom}f\cap \limfunc{dom}h%
\text{ and }s\geq 0, \\
+\infty , & \text{else.}%
\end{array}%
\right.
\end{equation*}%
Note that $\Delta _{1}\subset \Delta .$ For each $x\in X$ we have
\begin{equation*}
\sup\limits_{s\geq 0}L_{1}\left( x,s\right) =f(x)+\delta _{E}\left( x\right)
,
\end{equation*}%
and
\begin{equation*}
\inf_{x\in X}\sup_{s\geq 0}L_{1}\left( x,s\right) =\inf \left( P\right) .
\end{equation*}

The corresponding Lagrangian dual problem, say \textit{sup-dual problem},
reads
\begin{equation*}
(D_{1})\text{ }\sup_{s\geq 0}\inf_{x\in \Delta _{1}}\left( f(x)+sh(x)\right)
.
\end{equation*}%
Let us introduce the \textit{sup-value function} $v_{1}:\mathbb{%
R\longrightarrow }\overline{\mathbb{R}}$ associated with $\left( P\right) $
via $L_{1},$ namely,
\begin{equation*}
v_{1}\left( r\right) :=\inf \left\{ f(x):h(x)\leq r\right\} ,\text{ }r\in
\mathbb{R},
\end{equation*}%
which is non-increasing and satisfies
\begin{equation}
\overline{v}_{1}\left( 0\right) =\lim\limits_{\varepsilon \downarrow
0}v_{1}\left( \varepsilon \right) =\lim\limits_{\varepsilon \downarrow
0}\inf \left\{ f\left( x\right) :f_{t}\left( x\right) \leq \varepsilon ,t\in
T\right\} .  \label{3.2}
\end{equation}

\begin{lemma}
\label{lemma3.1} $\sup (D)\leq \sup (D_{1})\leq \inf \left( P\right) .$
\end{lemma}

\noindent \textbf{Proof} Let us prove the first inequality (the second being
obvious). Given $\lambda \in \mathbb{R}_{+}^{(T)},$ one has to check that
\begin{equation*}
\inf_{x\in \Delta }\left( f\left( x\right) +\sum_{t\in T}\lambda
_{t}f_{t}\left( x\right) \right) \leq \sup (D_{1}).
\end{equation*}

If $\limfunc{supp}\lambda =\emptyset ,$ then
\begin{equation*}
\inf_{x\in \Delta }\left( f\left( x\right) +\sum_{t\in T}\lambda
_{t}f_{t}\left( x\right) \right) =\inf\limits_{x\in \Delta }f\leq
\inf\limits_{x\in \Delta _{1}}f\leq \sup (D_{1})
\end{equation*}%
and we are done.

If $\limfunc{supp}\lambda \neq \emptyset ,$ one has, for $s=\sum_{t\in
T}\lambda _{t},$%
\begin{equation*}
\begin{array}{ll}
\sup (D_{1}) & \geq \inf_{x\in \Delta _{1}}\left( f(x)+sh(x)\right) \medskip
\\
& \geq \inf_{x\in \Delta _{1}}\left( f(x)+s\sum_{t\in T}\frac{\lambda _{t}}{s%
}f_{t}\left( x\right) \right) \medskip \\
& \geq \inf_{x\in \Delta _{1}}\left( f\left( x\right) +\sum_{t\in T}\lambda
_{t}f_{t}\left( x\right) \right) \medskip \\
& \geq \inf_{x\in \Delta }\left( f\left( x\right) +\sum_{t\in T}\lambda
_{t}f_{t}\left( x\right) \right) .%
\end{array}%
\end{equation*}%
\hfill \hfill $\square $

\begin{proposition}[Limiting formula for $\sup (D_{1})$]
\label{Prop3.1} Assume either $\overline{v}_{1}\left( 0\right) \neq +\infty $
or $\sup (D_{1})\neq -\infty .$ Then we have%
\begin{equation*}
\sup (D_{1})=\lim\limits_{\varepsilon \downarrow 0}\inf \left\{ f\left(
x\right) :f_{t}\left( x\right) \leq \varepsilon ,t\in T\right\} .
\end{equation*}
\end{proposition}

\noindent \textbf{Proof} By (\ref{3.2}), the right-hand side of (\ref{3.2})
coincides with $\overline{v}_{1}\left( 0\right) .$ By definition of $v_{1}$
we have (as for $v$), $v_{1}^{\ast \ast }\left( 0\right) =\sup (D_{1}).$
Since $v_{1}$ is convex and either $\overline{v}_{1}\left( 0\right) \neq
+\infty $ or $v_{1}^{\ast \ast }\left( 0\right) \neq -\infty ,$ we then
have, by \cite[Proposition 1]{Borwein80}, $\sup (D_{1})=\overline{v}%
_{1}\left( 0\right) $ and we are done.\hfill $\square \medskip $

\begin{proposition}[Limiting formula for $\inf \left( P\right) $]
\label{Prop3.2} Assume the strong Slater condition%
\begin{equation}
\exists \alpha >0,\text{ }\exists a\in \func{dom}f:\ \ f_{t}(a)\leq -\alpha
,\ \forall t\in T,  \label{3.3}
\end{equation}%
holds. Then we have%
\begin{equation}
\inf \left( P\right) =\max_{s\geq 0}\inf_{x\in \Delta _{1}}\left(
f(x)+sh(x)\right) =\lim\limits_{\varepsilon \downarrow 0}\inf \left\{
f\left( x\right) :f_{t}\left( x\right) \leq \varepsilon ,t\in T\right\} .
\label{3.4}
\end{equation}
\end{proposition}

\textbf{Proof} By definition of $h$ we have%
\begin{equation*}
\inf \left( P\right) =\inf \left\{ f(x):h(x)\leq 0\right\} .
\end{equation*}

Note that (\ref{3.3}) amounts to the usual Slater condition relative to $h:$%
\begin{equation*}
\exists a\in \func{dom}f:\ \ h(a)<0.
\end{equation*}

Since the functions $\ f$ and $h$ are convex, we then have (see, e.g., \cite[%
Lemma 1]{LV98})
\begin{equation*}
\inf \left( P\right) =\max_{s\geq 0}\inf_{x\in \Delta _{1}}\left(
f(x)+sh(x)\right) =\max \left( D_{1}\right) .
\end{equation*}

By (\ref{3.3}) we have $\overline{v}_{1}\left( 0\right) \leq v_{1}\left(
0\right) <+\infty .$ By Proposition \ref{Prop3.1} it follows that
\begin{equation*}
\sup (D_{1})=\lim\limits_{\varepsilon \downarrow 0}\inf \left\{ f\left(
x\right) :f_{t}\left( x\right) \leq \varepsilon ,t\in T\right\}
\end{equation*}%
and we are done.\hfill $\square \medskip $

Let us revisit Example \ref{Exam2}, where (\ref{3.4}) fails. Any candidate $%
a $ to be strong Slater point is feasible. Let $a$ be a feasible solution of
$\left( P\right) .$ Then $a=\left( a_{1},0\right) ,$ with $a_{1}\leq 0,$ and
$h\left( a\right) \geq \sup \left\{ t^{-1}a_{1}:t=3,4,...\right\} =0.$ Thus,
$h\left( a\right) =0$ and the strong Slater constraint qualification (\ref%
{3.3}) fails. However, by Proposition \pageref{Prop3.1}, we have%
\begin{equation*}
\sup (D_{1})=\lim\limits_{\varepsilon \downarrow 0}\inf \left\{ f\left(
x\right) :h\left( x\right) \leq \varepsilon \right\}
=\lim\limits_{\varepsilon \downarrow 0}-\varepsilon =0
\end{equation*}%
and, finally,%
\begin{equation*}
\begin{array}{c}
-1=\sup (D_{0})=\sup (D)<\sup (D_{1})=0=\min \left( P\right) \\
=\inf \left\{ f\left( x\right) :h\left( x\right) =0\right\}
=\lim\limits_{\varepsilon \downarrow 0}\inf \left\{ f\left( x\right)
:h\left( x\right) \leq \varepsilon \right\} .%
\end{array}%
\end{equation*}

\begin{remark}
\label{Rem 3.1} In the case when $T$ is finite, condition (\ref{3.3}) reads%
\begin{equation*}
\exists a\in \func{dom}f:\ \ f_{t}(a)<0,\ \forall t\in T,
\end{equation*}%
that is the familiar Slater constraint qualification. One has also $\Delta
_{1}=\left( \tbigcap\limits_{t\in T}\ \limfunc{dom}f_{t}\right) \cap
\limfunc{dom}f$ and, by Proposition \ref{Prop3.2}, there exists $\overline{s}%
\geq 0$ such that%
\begin{equation*}
\inf \left( P\right) =\inf_{x\in \Delta _{1}}\left( f(x)+\overline{s}%
h(x)\right) =\inf_{x\in \Delta _{1}}\sup\limits_{\nu \in S_{T}}\left( f(x)+%
\overline{s}\sum_{t\in T}\nu _{t}f_{t}\left( x\right) \right) ,
\end{equation*}%
where $S_{T}=\left\{ \nu \in \mathbb{R}_{+}^{T}:\sum_{t\in T}\nu
_{t}=1\right\} $ is the unit simplex in $\mathbb{R}^{T}.$ By the minimax
theorem \cite[Theorem 2.10.1]{Za02}, with $A=S_{T}$ and $B=\Delta _{1},$
there exists $\overline{\nu }\in S_{T}$ such that
\begin{equation*}
\inf \left( P\right) =\inf_{x\in \Delta _{1}}\left( f(x)+\overline{s}%
\sum_{t\in T}\overline{\nu }_{t}f_{t}\left( x\right) \right) \leq \sup
(D)\leq \inf \left( P\right)
\end{equation*}%
and, consequently, $\inf \left( P\right) =\max \left( D\right) ,$ which is
the strong duality theorem \cite[Theorem 2.9.3]{Za02} without assuming a
topological structure on the basic linear space $X$ (see also \cite[Remark 8]%
{LV21}).
\end{remark}

Concerning Example \ref{Exam1}, let us note that
\begin{equation*}
\max (D_{0})=0<1=\max (D)=\lim\limits_{\varepsilon \downarrow 0}\inf \left\{
f\left( x\right) :f_{1}\left( x\right) \leq \varepsilon \right\} =\min
\left( P\right) ,
\end{equation*}%
which also contradicts \cite[Proposition 3.1]{Karney83}.\medskip

\textbf{Acknowledgement }This research was partially supported by Ministerio
de Ciencia, Innovaci\'{o}n y Universidades (MCIU), Agencia Estatal de
Investigaci\'{o}n (AEI), and European Regional Development Fund (ERDF),
Project PGC2018-097960-B-C22.\medskip\


\begin{thebibliography}{99}
\bibitem{BDJ82} Blair, C.E., Duffin, R.J., Jeroslow, R.G.: A limiting
infisup theorem. J Optim Theory Appl \textbf{37}, 163-175 (1982)

\bibitem{Borwein80} Borwein, J.M.: A note on perfect duality and limiting
Lagrangeans. Math. Programming \textbf{18}, 330-337 (1980)

\bibitem{Bot10} Bo\c{t}, R.I.: Conjugate Duality in Convex Optimization.
Springer-Verlag, Berlin/Heidelberg (2010)

\bibitem{DGLV21} Dinh, N., Goberna, M.A., L\'{o}pez, M.A.: Relaxed
Lagrangian duality in convex infinite optimization: reverse strong duality
and optimality. Preprint. Available at http://arxiv.org/abs/2106.09299

\bibitem{Duffin73} Duffin, R.J.: Convex analysis treated by linear
programming. Math. Programming \textbf{4}, 125-143 (1973)

\bibitem{DJ79} Duffin, R.J., Jeroslow, R.G.: The Limiting Lagrangian.
Georgia Institute of Technology, Management Science Technical Reports No.
MS-79-13 (1979)

\bibitem{GL98} Goberna, M.A., L\'{o}pez, M.A.: Linear Semi-Infinite
Optimization. J. Wiley, Chichester, U.K., (1998)

\bibitem{Karney83} Karney, D.F.: A duality theorem for semi-infinite convex
programs and their finite subprograms. Math. Programming \textbf{27}, 75-82
(1983)

\bibitem{KM86} Karney, D.F., Morley, T.D.: Limiting Lagrangians: A primal
approach. J Optim. Theory Appl. 48, 163-174 (1986).

\bibitem{LV98} Lemaire, B., Volle, M.: Duality in DC programming.
Generalized convexity, generalized monotonicity: recent results (Luminy,
1996), 331-345, Nonconvex Optim. Appl. \textbf{27}, Kluwer, Dordrecht (1998)

\bibitem{LV21} Luc, D.T., Volle, M.: Algebraic approach to duality in
optimization and applications. Set-Valued Var. Anal., to appear.

\bibitem{Pomerol82} Pomerol, J.-Ch.: A note on limiting infisup theorems.
Math. Programming \textbf{30}, 238-241 (1984)

\bibitem{Rockafellar74} Rockafellar, R.T.: Conjugate Duality and
Optimization. SIAM, Philadelphia, P.A. (1974)

\bibitem{Za02} Z\u{a}linescu, C.: Convex analysis in general vector spaces.
World Scientific, River Edge, N.J. (2002)
\end{thebibliography}
\end{document}